\documentclass[12pt,twoside,reqno]{amsart}

\usepackage[OT1]{fontenc}    %
\usepackage{type1cm}         %
\usepackage{amsthm}          %
\usepackage[dvips]{graphicx} 

\usepackage{amsfonts}
\usepackage{amsmath}
\usepackage{amssymb}
\usepackage{wasysym}

\usepackage{verbatim}        

\numberwithin{equation}{section}

\theoremstyle{plain}
\newtheorem{theorem}{Theorem}[section]

\newtheorem{lemma}[theorem]{Lemma}

\theoremstyle{remark}

\newtheorem*{ack}{Acknowledgement}

\theoremstyle{definition}

\newcommand{\leb}{\mathcal{L}}
\newcommand{\supp}{\textrm{Supp}}

\newcommand{\e}{\varepsilon}
\newcommand{\D}{\mathcal{D}}

\newcommand{\I}{\mathcal{I}}

\newcommand{\C}{\mathcal{C}}

\newcommand{\T}{\mathcal{T}}
\newcommand{\wI}{\widetilde{I}}
\newcommand{\ophi}{\overline{\phi}}
\newcommand{\rot}{\mathbf{R}}
\newcommand{\RR}{\mathbb{R}}
\newcommand{\mh}{\mathcal{H}}
\newcommand{\XX}{\mathbb{X}}
\newcommand{\YY}{\mathbb{Y}}
\newcommand{\PP}{\mathbf{P}}

\DeclareMathOperator{\dist}{dist}
 
\title[Convolutions of Cantor measures]
{Convolutions of Cantor measures\\ without resonance}
\author{Fedor Nazarov, Yuval Peres and Pablo Shmerkin}

\begin{document}

\begin{abstract}
Denote by  $\mu_a$  the distribution of the random sum
$ \; \; (1-a) \sum_{j=0}^\infty \omega_j a^j$,
 where $\mathbf{P}(\omega_j=0)=\mathbf{P}(\omega_j=1)=1/2$ and all
the choices are independent. For $0<a<1/2$, the measure  $\mu_a$ is
  supported on $C_a$, the central Cantor set
obtained by starting with the closed united interval, removing an
open central interval of length $(1-2a)$, and iterating this
process inductively on each of the remaining intervals.
We investigate the convolutions $\mu_a * (\mu_b \circ
S_\lambda^{-1})$, where $S_\lambda(x)=\lambda x$ is a rescaling
map. We prove that if the ratio $\log b/\log a$ is irrational and
$\lambda\neq 0$, then
\[
D(\mu_a *(\mu_b\circ S_\lambda^{-1})) =
\min(\dim_H(C_a)+\dim_H(C_b),1),
\]
where $D$ denotes any of correlation, Hausdorff or packing
dimension of a measure.

We also show that, perhaps surprisingly, for uncountably many values of $\lambda$ the convolution $\mu_{1/4} *(\mu_{1/3}\circ S_\lambda^{-1})$ is a singular measure, although $\dim_H(C_{1/4})+\dim_H(C_{1/3})>1$ and $\log (1/3) /\log (1/4)$ is irrational.
\end{abstract}

\date{\today}
\thanks{P.S. acknowledges support from the Academy of Finland, Microsoft Research, EPSRC grant EP/E050441/1 and the University of Manchester.}

\maketitle

\section{Introduction and statement of results}

Given $0<a<1/2$, let $C_a$ be the Cantor set obtained by starting with the closed
unit interval, removing a central open interval of length $1-2a$,
and continuing this process inductively on each of the remaining
intervals. Formally,
\[
C_a = \left\{ (1-a) \sum_{j=0}^\infty \omega_j a^j :
\omega_j\in\{0,1\} \textrm{ for all } j \right\}.
\]

The set $C_a$ supports a natural probability measure $\mu_a$ which
assigns mass $2^{-n}$ to each interval of length $a^n$ in the
construction. The measure $\mu_a$ can be defined in several
alternative ways. For example, it is the normalized restriction of
$\dim_H(C_a)$-dimensional Hausdorff measure to $C_a$. It is also
the distribution of the random infinite sum
\begin{equation} \label{eq:defrandomsum}
(1-a) \sum_{j=0}^\infty \omega_j a^j,
\end{equation}
where $\mathbf{P}(\omega_j=0)=\mathbf{P}(\omega_j)=1/2$ and all
choices are independent. These equivalences are well known and
easy to verify.

In this paper we study convolutions of the form $\mu_a *
(\mu_b\circ S_\lambda^{-1})$, where $S_\lambda(x)=\lambda x$
scales by a factor of $\lambda$. We will show that, under a
natural irrationality condition, $\mu_a * (\mu_b\circ
S_\lambda^{-1})$ has ``fractal dimension'' equal to the sum of the
Hausdorff dimensions of $C_a$ and $C_b$, provided this is at most
one. What we mean for fractal dimension is made precise below; we
will in fact show that this is true for several commonly used
concepts of dimension of a measure.

The study of these convolutions goes back to Senge and Straus
\cite{SengeStraus1971} who, answering a question that Salem posed
in \cite{Salem1963}, characterized all the pairs $a,b$ such that
$\phi_{a,b}(\xi)\nrightarrow 0$ as $\xi\rightarrow\infty$, where
$\phi_{a,b}$ denotes the Fourier transform of $\mu_a * \mu_b$.
Senge and Straus showed that this happens only if $1/a$ and $1/b$
are Pisot numbers and $\log b/\log a$ is rational (Recall that a
Pisot number is an algebraic integer larger than one, such that
all its algebraic conjugates are smaller than one in modulus).

Let us write
\[
\nu_{a,b}^\lambda = \mu_a * (\mu_b\circ S_\lambda^{-1}).
\]
There are two sharply different cases in the study of
$\nu_{a,b}^\lambda$: The subcritical case $d_a+d_b < 1$ and the supercritical case $d_a+d_b>1$ (the critical case $d_a+d_b=1$ is often analyzed separately).

In the subcritical case, the measure $\nu_{a,b}^\lambda$ is always
singular, as it is supported on $C_a+ \lambda C_b$, which has
Hausdorff dimension at most $d_a+d_b<1$ (see
\eqref{eq:ineqs_sumdims} below). Thus, in this case the interest
lies in the degree of singularity of $\nu_{a,b}^\lambda$, as
measured by some concept of fractal dimension. In particular, one
is interested in whether there is a ``dimension drop'', i.e.
whether the dimension of $\nu_{a,b}^\lambda$ is strictly smaller
than the dimension of $\mu_a \times \mu_b$. We will prove that if
$\log b/\log a$ is irrational, then there is no dimension drop for
any $\lambda\neq 0$, for several different concepts of dimension
of a measure; see Theorem \ref{th:main_result} and the discussion
afterwards. One motivation comes from the results in
\cite{PeresShmerkin2009}, where it is proved that for all pairs
$0<a,b<1/2$ such that $\log b/\log a$ is irrational and all
$\lambda\neq 0$,
\[
\dim_H(C_a+\lambda C_b) = \min(\dim_H(C_a)+\dim_H(C_b),1),
\]
where $\dim_H$ stands for Hausdorff dimension. The proofs in
\cite{PeresShmerkin2009} involve the construction of an ad-hoc
measure supported on $C_a+\lambda C_b$, which is not related in a natural
way to $\nu_{a,b}^\lambda$. In this paper we base the arguments on more conceptual ergodic-theoretical ideas.

In the case $d_a+d_b>1$,  one would expect $\nu_{a,b}^\lambda$ to be absolutely
continuous as long as $\log b/\log a\notin \mathbb{Q}$. However,
we will show that this is not always the case. More precisely, we
will prove that whenever $1/a$ and $1/b$ are Pisot numbers and $\log b/\log a$ is irrational, there is a
dense $G_\delta$ set of parameters $\lambda$, such that the
Fourier transform of $\nu_{a,b}^\lambda$ does not go to zero at
infinity; see Theorem \ref{th:example} in Section \ref{sec:largedimension}.

In order to state our main result about the fractal dimensions of
$\nu_{a,b}^\lambda$, we start by recalling the definition of
correlation dimension of a measure. Given a Borel measure $\nu$ on
$\mathbb{R}^n$ and $r>0$, let
\[
\mathcal{C}_\nu(r) = \int \nu (B(x,r)) d\nu(x),
\]
where $B(x,r)$ denotes the closed ball with center $x$ and radius
$r$. The \textbf{lower correlation dimension} of $\nu$ is defined
as
\[
\underline{D}(\nu) = \liminf_{r\downarrow 0} \frac{\log
\mathcal{C}_\nu(r)}{\log r}.
\]
The \textbf{upper correlation dimension} $\overline{D}(\nu)$ is
defined analogously by taking the $\limsup$. If
$\underline{D}(\nu) = \overline{D}(\nu)$ we say that the
\textbf{correlation dimension} $D(\nu)$ exists, and is given by
the common value.  Other definitions of correlation dimension are
often used. For example, the lower correlation dimension of $\nu$
is the supremum of all $\alpha$ such that
\begin{equation} \label{eq:def_energy}
I_\alpha(\nu) = \int\int |x-y|^{-\alpha} d\nu(x)d\nu(y) < \infty.
\end{equation}
This is well-known, see e.g. \cite[Proposition 2.3]{SauerYorke1997}.
We can now state our main result:

\begin{theorem} \label{th:main_result}
Let $0<a,b<1$. If $\log b/\log a$ is irrational, then
\begin{equation} \label{eq:main_result}
D(\nu_{a,b}^\lambda) = \min(d_a+d_b,1),
\end{equation}
for all $\lambda \neq 0$.
\end{theorem}

Let us make some comments on this statement. Observe that the
measure $\nu_{a,b}^\lambda$ is, up to affine equivalence, the
push-down of the product measure $\mu_a\times \mu_b$ by orthogonal
projection onto the line
\[
\{ t(\cos(\theta),\sin(\theta)):t\in\mathbb{R}\},
\]
where $\theta = \arctan(\lambda)$. The potential-theoretic proof
of Marstrand's projection theorem (see e.g. \cite[Chapter
9]{Mattila1995}) implies that \eqref{eq:main_result} holds for
almost every $\lambda$. Our contribution is to prove this for
\textit{every} $\lambda \neq 0$, and in particular for
$\lambda=1$, which yields the correlation dimension of the
convolution measure $\mu_a * \mu_b$.

Correlation dimension is just one of the several concepts of
dimension of a measure which are often used. For example, the
Hausdorff and packing dimensions of a finite Borel measure $\nu$
on $\mathbb{R}^n$ are defined as
\begin{eqnarray*}
\dim_H(\nu) & = & \inf\{ \dim_H(E) : \nu(\mathbb{R}^n\backslash
E)=0 \},\\
\dim_P(\nu) & = & \inf\{ \dim_P(E) : \nu(\mathbb{R}^n\backslash
E)=0 \}.
\end{eqnarray*}
(Here $\dim_P$ denotes packing dimension; see e.g. \cite[Chapter
5]{Mattila1995} for its definition and basic properties). It is
then an easy consequence of Frostman's lemma (see \cite[Chapter
8]{Mattila1995}) and the definitions that
\begin{equation} \label{eq:ineqs_hausdorff}
\underline{D}(\nu)  \le \dim_H(\nu) \le \dim_H(\supp(\nu)) \le n,
\end{equation}
where $\supp$ denotes the support of the measure. It is also
immediate that
\begin{equation} \label{eq:ineqs_packing}
\dim_H(\nu) \le \dim_P(\nu) \le \dim_P(\supp(\nu)) \le n.
\end{equation}
Since orthogonal projections do not increase either Hausdorff or
packing dimension, and $C_a\times C_b$ has Hausdorff and packing
dimension $d_a+d_b$, we get
\begin{equation} \label{eq:ineqs_sumdims}
\dim_H(C_a+\lambda C_b) \le \dim_P(C_a+\lambda C_b) \le
\min(d_a+d_b,1)
\end{equation}
for \textit{all} $0<a,b<1/2$. Hence we deduce from Theorem
\ref{th:main_result}, \eqref{eq:ineqs_hausdorff} and
\eqref{eq:ineqs_packing} that, whenever $\log b/\log a$ is
irrational,
\begin{equation} \label{eq:dimsconvmeasure}
D(\nu_{a,b}^\lambda) = \dim_H(\nu_{a,b}^\lambda) =
\dim_P(\nu_{a,b}^\lambda) = \min(d_a+d_b,1).
\end{equation}

Given a measure $\nu$ on Euclidean space,
the \textbf{local dimensions} of $\nu$ are defined as
\[
\dim_{\textrm{loc}}(\nu)(x) = \lim_{r\downarrow 0} \frac{\log(\nu(B(x,r)))}{\log(r)},
\]
whenever the limit exists; otherwise one speaks of lower and upper
local dimensions. When the local dimension exists and is constant
$\nu$-almost everywhere, it is said that $\nu$ is \textbf{exact
dimensional}. Always assuming that $\log b/\log a$ is irrational,
it follows from \eqref{eq:dimsconvmeasure} and \cite[Theorems 1.2
and 1.4]{FanLauRao2002} that
\[
\dim_{\textrm{loc}}(\mu_a*\mu_b)(x) = \min(d_a+d_b,1) \quad \textrm{ for } (\mu_a*\mu_b) \textrm{-a.e.} x.
\]
We remark that in general all of the concepts of dimension of a
measure that we discussed can differ; even if a measure is
exact-dimensional, it may happen that its correlation dimension is
strictly smaller than the almost-sure value of the local
dimension.

We note that if $\log b/\log a$ is rational, and $d_a+d_b\le 1$, then
\[
\dim_H(\mu_a*\mu_b) \le \dim_H(C_a+C_b) < d_a+d_b.
\]
(See \cite{PeresShmerkin2009} for a proof). Borrowing the
terminology of \cite{PeresShmerkin2009}, we can summarize our
results as saying that \textbf{algebraic resonance}, defined as
the rationality of $\log b/\log a$, is equivalent to
\textbf{geometric resonance} for the measures $\mu_a$ and $\mu_b$,
defined by the condition that there is a dimension drop for at least one orthogonal projection in a non-principal direction.

Compared to \cite{PeresShmerkin2009}, one basic new ingredient in our proofs is the
construction of a subadditive cocycle reflecting the structure of
the orthogonal projections of $\mu_a\times\mu_b$ at different
scales and for different angles. This is obtained by adapting the
proof of the existence of $L^q$ dimensions for self-conformal
measures in \cite{PeresSolomyak2000}. We will also make use of a theorem of A. Furman on subadditive cocycles over a uniquely ergodic transformation.

The paper is structured as follows. In Section \ref{sec:notation}
we introduce notation and prove several lemmas in preparation for
the proof of Theorem \ref{th:main_result}, which is contained in
Section \ref{sec:proof}. In Section \ref{sec:largedimension}, we
prove the singularity of $\nu_{1/3,1/4}^\lambda$
for an uncountable set of $\lambda$, and discuss some
implications. We finish the paper with some remarks and open
questions in Section \ref{sec:conclusion}.

\section{Notation and preliminary lemmas}
\label{sec:notation}

In this section we collect several definitions and lemmas which
will be used in the proof of Theorem \ref{th:main_result}.

\subsection{Notation and basic facts}

From now on we will assume that $0<a<b<1/2$ and $\log b/\log a$ is
irrational.

Consider the product set $E =C_a\times C_b$, and let $\eta =
\mu_a\times\mu_b$. This is, up to a constant multiple,
$(d_a+d_b)$-dimensional Hausdorff measure restricted to $E$. It is
well known, and easy to verify, that $\eta$ is Ahlfors-regular,
i.e.
\begin{equation} \label{eq:ahlfors_regular}
C^{-1} r^{d_a+d_b} \le \eta(B(x,r)) \le C r^{d_a+d_b},
\end{equation}
for some constant $C>1$, all $r>0$ and all $x\in E$. This implies
that
\begin{equation} \label{eq:dimension_eta}
D(\eta) = \dim_H(\eta) = \dim_P(\eta) = d_a+d_b.
\end{equation}

The set $C_a$ can also be realized as the attractor of the
iterated function system $\{ f_{a,0}, f_{a,1}\}$, where
\begin{equation} \label{eq:ifs}
f_{a,i}(x) = a x + i(1-a).
\end{equation}
In other words,
\[
C_a = f_{a,0}(C_a) \cup f_{a,1}(C_a).
\]
Likewise, the measure $\mu_a$ satisfies the self-similarity
relation
\begin{equation} \label{eq:selfsimilarity}
\mu_a(F) =
\frac{1}{2}\mu_a(f_{a,0}^{-1}(F))+\frac{1}{2}\mu_a(f_{a,1}^{-1}(F)),
\end{equation}
where $F$ is an arbitrary Borel set. This is well known; when $F$
is a basic interval in the construction of $C_a$ it follows from
the scaling property of Hausdorff measure and self-similarity, and
the case of general $F$ is obtained by noting that basic intervals
in the construction of $C_a$ form a basis of closed subsets of
$C_a$.

If $u=(u_1,\ldots,u_k)\in\{0,1\}^k$, let
\[
f_{a,u} = f_{a,u_1}\circ \cdots\circ f_{a,u_k};
\]
analogously one defines $f_{b,u}$. Let
\[
\XX = \bigcup_{k,l=0}^\infty \{0,1\}^k\times \{0,1\}^l.
\]
Given $\xi = (u,v)\in \XX$ we will denote
\[
f_\xi(x,y) = (f_{a,u}(x),f_{b,v}(y)).
\]
Moreover, we will write $E(\xi)=f_\xi(E)$ and $Q(\xi) = f_\xi(Q)$,
where $Q=I\times I$ is the unit square. If $\xi_i = (u_i,v_i)$,
$i=1,2$, by $\xi_1 \xi_2$ we will mean the juxtaposition $(u_1
u_2, v_1 v_2)$. The empty word will be denoted by $\varnothing$.
Finally, if $\xi=(u,v)$, the \textbf{length pair} $|\xi|$ of $\xi$
is the pair $(|u|,|v|)$, where $|u|,|v|$ are the lengths of the
corresponding words.

Let $\ell$ be any integer such that $b/a < 1/b^\ell$. Write
$\alpha = \log(b/a)$, $\beta = \log(1/b^\ell)$; notice that
$0<\alpha<\beta$, and that $\alpha/\beta$ is irrational because of
our assumption that $\log b /\log a$ is irrational. Moreover, note
that $\beta$ can be made arbitrarily large by starting with an
appropriately large $\ell$. Endow $[0,\beta)$ with normalized
Lebesgue measure $\leb$. Let also $\rot:[0,\beta)\rightarrow
[0,\beta)$ be given by
\[
\rot(x) = x + \alpha \mod(\beta).
\]
The irrationality of $\alpha/\beta$ implies that this is a
uniquely ergodic transformation. This fact will be crucial in the
proof: it is the only place where the irrationality of $\log
b/\log a$ is used.

We inductively construct two families $\{ \XX_n \}_{n=0}^\infty$,
$\{ \YY_n \}_{n=0}^\infty$ of subsets of $\XX$. We set $ \XX_0 =\{
(\varnothing,\varnothing)\}$ (recall that $\varnothing$ denotes
the empty word). Once $ \XX_n$ has been defined, we define $
\XX_{n+1}$ in the following way:
\begin{itemize}
\item if $\rot^n(0) + \alpha < \beta$, then
\[
 \XX_{n+1} = \{ (\xi) (i,k) : \xi\in \XX_n,\, i,k\in \{0,1\} \}.
\]
\item If $\rot^n(0) + \alpha > \beta$, then
\[
 \XX_{n+1} = \{  (\xi) (i,v) : \xi\in \XX_n,
\,i\in\{0,1\},v\in\{0,1\}^{\ell+1} \}.
\]
\end{itemize}
Further, we let
\[
\YY_n = \{ (\xi)(\varnothing,v) : \xi \in \XX_n, v\in\{0,1\}^\ell
\}.
\]
One can readily check, from this inductive construction and the
definition of $\rot$, that the following properties are satisfied:
\begin{enumerate}
\item[(I)] If $\xi\in \XX_n$, then $Q(\xi)$ is a rectangle  of size
\[
a^n \times  a^n \exp(\rot^n(0)).
\]
If $\xi'\in \YY_n$, then $Q(\xi')$ is a rectangle of size
\[
a^n \times a^n \exp(\rot^n(0)-\beta).
\]
In particular, all elements of $\XX_n$ have the same length pair,
as do all elements of $\YY_n$.
\item[(II)] The cylinders based at elements of $\XX_n$ form a partition of the symbolic space. In other words, if
$\omega,\omega'\in\{0,1\}^\mathbb{N}$ are infinite sequences, then
there is exactly one $\xi = (u,v)\in \XX_n$ such that $\omega$
starts with $u$ and $\omega'$ starts with $v$. The same holds for
$\YY_n$.
\end{enumerate}

The rectangles $Q(\xi)$ are cartesian products of basic intervals
of $C_a$ and $C_b$; by the first property, the logarithm of the ratio of the
lengths stays bounded and, moreover, behaves like an irrational
rotation. Property (II) guarantees that $\{ Q(\xi)\in \XX_n\}$ is
an efficient covering of $E$, and likewise for $\YY_n$.

Heuristically, for $\XX_n$, we start from the unit square, and
then at each inductive step we always go one level further in the
construction of $C_a$. With respect to the construction of $C_b$,
we go one level further for as long as the eccentricity of the
rectangles $Q(\xi)$ stays below $e^\beta=b^{-\ell}$ (note that,
since $b>a$, going one step further in both constructions has the
effect of increasing the eccentricity); otherwise, we go $\ell+1$
levels further, which has the effect of reducing the eccentricity
of $Q(\xi)$ back to a value between $1$ and $e^\beta$. For the
construction of $\YY_n$, we start from $\XX_n$ and go $\ell$
levels further in the construction of $C_b$, while keeping the
same basic intervals in the construction of $C_a$; this yields
rectangles $Q(\xi')$ with width greater than height.

Let $h_s(x,y) = (x,e^s y)$, and $\Pi(x,y)=x+y$. We will write
$\Pi_s = \Pi  h_s$. The assignment $s\rightarrow h_s$ is an action
of the additive group of real numbers by linear bijections of
$\RR^2$.

The following observation will prove very useful: let $\xi\in
\XX_n$. Then $f_\xi$ can be decomposed as
\begin{equation} \label{eq:decomposition_f}
f_\xi(x) =  a^n h_{\rot^n(0)}(x) + d_\xi,
\end{equation}
where $d_\xi$ is a translation vector in $\RR^2$. From this we
also obtain
\begin{equation} \label{eq:decomposition_f_inverse}
 f_\xi^{-1}(x) = a^{-n} h_{-\rot^n(0)}(x) + \hat{d}_\xi,
\end{equation}
where $\hat{d}_\xi = - a^{-n} h_{-\rot^n(0)} (d_\xi)$. Of course,
similar decompositions are valid for $f_{\xi'}$, $f_{\xi'}^{-1}$
for $\xi'\in \YY_n$.

We will denote the projected measure $\eta\circ\Pi_s^{-1}$ by
$\eta_s$. Notice that $\eta_s = \mu_a * (\mu_b \circ
S_{e^s}^{-1})$, where $S_\lambda(x) = \lambda x$. In particular,
$\eta_0 = \mu_a
* \mu_b$.

For $n\in\mathbb{N}$ we will denote
\[
\D_n = \{ [ j a^n, (j+1) a^n) :j\in\mathbb{Z} \}.
\]
Further, for $s\in\mathbb{R}$, we will denote by $\D_n(s)$ the
subset of $\D_n$ comprising all intervals which intersect
$\Pi_s(E)$.

Given $s\in\RR$, let us define
\begin{equation} \label{deftau}
\tau_n(s) = \sum_{I\in\D_n} \eta_s(I)^2 = \sum_{I\in \D_n(s)}
\eta_s(I)^2.
\end{equation}
These functions will be a useful discrete analogue of
$\mathcal{C}_{\eta_s}(a^n)$; indeed, both quantities are
comparable up to a multiplicative constant.

The next definition is adapted from \cite{PeresSolomyak2000}. Given $n \in\mathbb{N}$ and
$s\in\mathbb{R}$, we will say that a family of disjoint intervals
$\C$ is $(n, s)$-\textbf{good} if it is a minimal covering of
$\Pi_s(E)$ (meaning that no proper subset is a covering of
$\Pi_s(E)$), and each interval has length $a^n$.

\begin{lemma} \label{lemma:equivcovering1}
If $\C$ is $(n, s)$-good, then
\[
4^{-1}\tau_n(s) \le \sum_{C\in\C} \eta_s(C)^2 \le 4 \tau_n(s).
\]
\end{lemma}
\begin{proof}
We only prove the right-hand inequality since the left-hand
inequality is analogous. Since $\C$ is $(n, s)$-good, each element
$C$ of $\C$ intersects at most $2$ elements of $\D_n(s)$, say
$D(C,i), i=1,2$. Since $\D_n(s)$ is a covering of $\Pi_s(E)$,
\[
C \cap \Pi_s(E) \subset \left( D(C,1) \cup D(C,2)\right) \cap
\Pi_s(E).
\]
Using this and Cauchy-Schwartz,
\[
\eta_s(C)^2 \le \left(\sum_{i=1}^{2} \eta_s(D(C,i))\right)^2 \le 2
\sum_{i=1}^2 \eta_s(D(C,i))^2.
\]
Therefore
\[
 \sum_{C\in\C} \eta_s(C)^2 \le 2 \sum_{C\in C} \sum_{i=1}^2 \eta_s(D(C,i))^2.
\]
But since each element of $\D_n$ intersects at most $2$ elements
of $\C$, any element of $\D_n$ appears at most twice on the
right-hand side, and the lemma follows.
\end{proof}

\begin{lemma} \label{lemma:equivcovering2}
For any $s\in \mathbb{R}$, $m,n\in\mathbb{N}$ and $\xi\in \XX_n,
\xi'\in \YY_n$, the families
\begin{eqnarray*}
\I & = & \left\{ \Pi_{\rot^n(0)+s} f_\xi^{-1} \Pi_s^{-1} I : I\in
\D_{m+n}, I\cap \Pi_s(E(\xi))\neq\varnothing \right\}, \\
\I' & = & \left\{ \Pi_{\rot^n(0)+s-\beta} f_{\xi'}^{-1} \Pi_s^{-1}
I : I\in \D_{m+n}, I\cap \Pi_s(E(\xi'))\neq\varnothing \right\},
\end{eqnarray*}
are $(m, \rot^n(0)+s)$-good and $(m,\rot^n(0)+s-\beta)$-good,
respectively.
\end{lemma}
\begin{proof}
We will prove only that $\I$ is $(m,\rot^n(0)+s)$-good, since the
proof for $\I'$ is analogous. Let $t=\rot^n(0)+s$. Observe that
\[
f_\xi(E) = E(\xi)\subset \bigcup_{I\in
\D_{m+n},I\cap\Pi_s(E(\xi))\neq\varnothing} \Pi_s^{-1}(I),
\]
and from this it follows that $\I$ is a covering of $\Pi_t(E)$.
Using \eqref{eq:decomposition_f_inverse}, linearity of $\Pi_t$,
and the action properties of $s\rightarrow h_s$, we get
\begin{eqnarray*}
\Pi_t f_\xi^{-1} \Pi_s^{-1} I & = & \Pi_t \hat{d}_\xi + a^{-n} \Pi_t h_{-\rot^n(0)} \Pi_s^{-1}(I)\\
 & = & \Pi_t \hat{d}_\xi + a^{-n} \Pi  h_t  h_{-\rot^n(0)}  h_s^{-1} \Pi^{-1}(I)\\
& = & \Pi_t \hat{d}_\xi + a^{-n}  I.
\end{eqnarray*}
This shows that the elements of $\I$ are pairwise disjoint (since
the elements of $\D_{m+n}$ are), and have length $a^m$. Since, by
definition, all elements of $\I$ intersect $\Pi_t(E)$, the
covering $\I$ is optimal. This completes the proof of the lemma.
\end{proof}

We finish this section by recalling the result of Furman alluded
to in the introduction.

\begin{theorem} \label{th:furman}
Let $X$ be a compact metric space, and let $T$ be a continuous
homeomorphism of $X$ with a unique invariant probability measure
$\mu$ (which must be ergodic). Further, let $\{ \phi_n\}$ be a
continuous subadditive cocycle over $(X,T)$. In other words, each
$\phi_n$ is a continuous real valued function on $X$, and
\[
\phi_{m+n}(x) \le \phi_m(x) + \phi_n(T^m x),
\]
for all $x\in X$. Let
\[
A = \lim_{n\rightarrow\infty} \frac{1}{n}\int_X \phi_n(x) d\mu.
\]
Note that $A$ is well defined by subadditivity. Then for almost
all $x\in X$,
\[
\lim_{n\rightarrow\infty} \frac{\phi_n(x)}{n} = A,
\]
and for all $x\in X$,
\[
\limsup_{n\rightarrow\infty} \frac{\phi_n(x)}{n} \le A.
\]
\end{theorem}
\begin{proof}
The first assertion is Kingman's subadditive ergodic theorem
\cite{Kingman1968}. The second one is proved in \cite[Theorem
1]{Furman} .
\end{proof}

\section{Proof of Theorem \ref{th:main_result}}
\label{sec:proof}

The proof consists of two main parts: in the first we construct a
submultiplicative cocycle over $([0,\beta),\rot)$ related to the
growth of $\tau_n(s)$. In the second part, we apply Theorem
\ref{th:furman} to this cocycle and deduce Theorem
\ref{th:main_result} from the potential-theoretic proof of
Marstrand's projection theorem.

Before starting the proof, we remark that due to the symmetry of
$C_a \times C_b$, we only need to show that \eqref{eq:main_result}
holds for all $\lambda \ge 1$ (which corresponds to orthogonal
projections for angles in $[\pi/4,\pi/2)$). Moreover, since
$\beta$ is arbitrarily large, it will be enough to establish
\eqref{eq:main_result} for all $\lambda \in [1,e^\beta)$.

\textbf{A submultiplicative cocycle}. The goal of this part is to
show that $\tau_n(\cdot)$ defined in (\ref{deftau}) satisfy
\begin{equation} \label{eq:submultiplicative}
\tau_{m+n}(s) \le A \tau_n(s) \tau_m(\rot^n(s)),
\end{equation}
for some $A>1$ independent of $m,n\in\mathbb{N}$ and $s\in
[0,\beta)$. In order to do this, we will follow the pattern of the
proof of existence of $L^q$ dimension in \cite{PeresSolomyak2000}
in the case $q>1$.

We will consider two cases, depending on whether
$\rot^n(0)+s<\beta$ or $\rot^n(0)+s > \beta$. In the first case,
we have $\rot^n(s) = \rot^n(0)+s$, while in the second, $\rot^n(s)
= \rot^n(0)+s-\beta$. In the proof of the first case we will use
the families $\{ \XX_n\}$, while the proof of the second is based
on the families $\{ \YY_n\}$. We will in fact only prove the first
case; the second follows in the same way, so details are left to
the reader.

Let $m,n\in\mathbb{N}$, and pick some $\xi\in  \XX_n$. Fix $s\in
[0,\beta)$, and let
\[
t=\rot^n(s) = \rot^n(0)+s.
\]
It follows from Lemmas \ref{lemma:equivcovering1} and
\ref{lemma:equivcovering2} that
\begin{equation} \label{eq:tauineq2}
\sum_{I\in\D_{m+n}, I\cap \Pi_s(E(\xi))\neq\varnothing}
\eta_t(\Pi_t f_\xi^{-1}\Pi_s^{-1} (I))^2 \le 4 \tau_m(t).
\end{equation}

We claim that
\begin{equation} \label{eq:claimstrip}
\Pi_t^{-1}\Pi_t f_\xi^{-1}\Pi_s^{-1}(I) = f_\xi^{-1}\Pi_s^{-1}(I).
\end{equation}
To see this, foliate $\mathbb{R}^2$ by fibers $\{\Pi_t^{-1}(x)\}$,
and note that $\Pi_t^{-1}\Pi_t(F)=F$ if and only $F$ contains the
fiber through all of its points. In the particular case
$F=f_\xi^{-1}\Pi_s^{-1}(I)$, we have:
\begin{eqnarray*}
p\in F & \Longleftrightarrow  & \Pi_s f_\xi(p) \in I \\
& \Longleftrightarrow &
\Pi h_s (a^n h_{R^n(0)}(p) + d_\xi) \in I \\
& \Longleftrightarrow & \Pi (a^n h_t(p)+ h_s d_\xi) \in I\\
& \Longleftrightarrow & a^n \Pi_t(p) + \Pi_s d_\xi \in I,
\end{eqnarray*}
where we used \eqref{eq:decomposition_f} in the second displayed
line, and the linearity of $\Pi$ in the fourth. Since the value of
$\Pi_t(p)$ is constant on the leaf through $p$,
\eqref{eq:claimstrip} is proved.

Notice that if $I\cap \Pi_s(E(\xi))=\varnothing$ then
$f_\xi^{-1}\Pi_s^{-1}(I)\cap E=\varnothing$. Combining this with
\eqref{eq:tauineq2} and \eqref{eq:claimstrip} yields
\begin{equation} \label{eq:tauineq3}
\sum_{I\in\D_{m+n}} \eta(f_\xi^{-1}\Pi_s^{-1} (I))^2 \le 4
\tau_m(t).
\end{equation}

Let $\D'_n$ be the family of unions of two consecutive elements of
$\D_n$; in other words,
\[
\D'_n = \{ [j a^n, (j+2) a^n) : j\in\mathbb{Z} \}.
\]
Let also $\D''_n = \D_n \cup \D'_n$. To each $I\in\D_{m+n}$ we
associate an element $\wI\in \D''_n$ in the following way: if $I$
is contained in an element of $\D_n$, let $\wI$ be this element;
otherwise, $I$ is contained in a unique element of $\D'_n$, and we
let $\wI$ be this element.

Iterating \eqref{eq:selfsimilarity} and using properties (I)-(II)
of $\{\XX_n\}$, we see that
\begin{equation} \label{eq:selfsimilarmeasure}
 \eta(F) = |\XX_n|^{-1} \sum_{\xi\in \XX_n}
 \eta(f_\xi^{-1} F),
\end{equation}
for any Borel set $F$. Using this we obtain
\begin{equation} \label{eq:measureiteration}
\eta_s(I) = \frac{1}{| \XX_n|} \sum_{\xi\in \XX_n}
\eta(f_\xi^{-1}\Pi_s^{-1} I) = \frac{1}{| \XX_n|} \sum_{\xi\in
\XX_n : \wI\cap \Pi_s(E(\xi))\neq\varnothing}
\eta(f_\xi^{-1}\Pi_s^{-1} I).
\end{equation}
For $J\in\D''_n$, write
\[
 \XX_n(J,s) = \{ \xi\in \XX_n : J\cap
\Pi_s(E(\xi))\neq\varnothing\}.
\]
From \eqref{eq:measureiteration}, an application of Cauchy-Schwartz yields
\[
\eta_s^2(I) \le | \XX_n|^{-2} | \XX_n(\wI,s)| \sum_{\xi\in
 \XX_n(\wI,s) } \eta(f_\xi^{-1} \Pi_s^{-1} I)^2.
\]
For a fixed $J\in \D''_n$, we add over all $I\in D_{m+n}$ such
that $\wI=J$, to get
\begin{eqnarray*}
\sum_{I\in\D_{m+n}, \wI=J} \eta^2_s(I)&  \le &  | \XX_n|^{-2} |
\XX_n(J,s)| \sum_{\xi\in  \XX_n(J,s) } \sum_{\wI=J}
\eta(f_\xi^{-1}
\Pi_s^{-1} I)^2 \\
& \le &  4\tau_m(t) | \XX_n|^{-2} | \XX_n(J,s)|^2,
\end{eqnarray*}
where in the last inequality we applied \eqref{eq:tauineq3}.
Adding over all $J\in\D''_n$, we obtain
\begin{equation} \label{eq:tauineq4}
\tau_{m+n}(s) \le 4 \tau_m(t) |\XX_n|^{-2} \sum_{J\in\D''_n} |
\XX_n(J,s)|^2.
\end{equation}

All the maps $\Pi_s$, $s\in [0,\beta)$, are Lipschitz with
Lipschitz constant uniformly bounded by $2 e^\beta$. Also, for $\xi\in \XX_n$, the diameter of $Q(\xi)$ is bounded by
$\sqrt{1+e^{2\beta}} a^n < 2 e^\beta a^n$. Given $J\in\D''_n$,
denote by $\hat{J}$ the interval with the same center as $J$
and length $|J|+16 e^{\beta} a^n$. By our previous observations,
if $\xi\in \XX_n(J,s)$, then $\Pi_s^{-1}(\hat{J})$ contains a ball
of radius $2e^\beta a^n$ centered at $E(\xi)$, and this implies
that $E \subset f_\xi^{-1}\Pi_s^{-1}(\hat{J})$. It follows that
\begin{eqnarray*}
| \XX_n|^{-1} | \XX_n(J,s)| & \le & \eta(E)^{-1} | \XX_n|^{-1}
\sum_{\xi\in \XX_n(J,s)}
\eta\left(f_\xi^{-1} \Pi_s^{-1} (\hat{J})\right)\\
&  \le & \eta(E)^{-1}\eta_s(\hat{J}) \le \eta(E)^{-1}\sum_{J'\in
\D_n, \dist(J',J)<8e^\beta a^n} \eta_s(J'),
\end{eqnarray*}
where we used \eqref{eq:selfsimilarmeasure}. Using Cauchy-Schwartz
once again, we deduce that
\[
\left(| \XX_n|^{-1} | \XX_n(J,s)|\right)^2 \le K_1 \sum_{J'\in
\D_n, \dist(J',J)< 8e^\beta a^n} \eta_s(J')^2,
\]
where $K_1 = \eta(E)^{-2}\lceil 16 e^\beta+2 \rceil$. Adding over
all $J \in\D''_n$, note that each element of $\D_n$ on the
right-hand side appears at most $2 K_1$ times, whence
\[
|\XX_n|^{-2} \sum_{J \in\D''_n} | \XX_n(J,s)|^2 \le 2 K_1^2
\tau_n(s).
\]
Together with \eqref{eq:tauineq4} this yields
\eqref{eq:submultiplicative}, as desired.

\medskip

\textbf{Conclusion of the proof}. Recall the definition of
$\mathcal{C}_\nu(r)$ given in the introduction. Let
\[
\phi_n(s) = \mathcal{C}_{\eta_s}(a^n).
\]
It is clear that the correlation dimension of $\eta_s$ exists if
and only if the limit $L=\lim_{n\rightarrow\infty}
\log\phi_n(s)/n$ exists, in which case $D(\eta_s)=L/\log(a)$.

Let us rewrite $\phi_n$ as
\begin{eqnarray*}
\phi_n(s) & = & \int \eta_s\left(\Pi_s(y)-a^n,\Pi_s(y)+a^n\right)
d\eta(y) \\
& = & \int \eta\left(\Pi_s^{-1}(\Pi_s(y)-a^n,\Pi_s(y)+a^n)\cap
Q\right) d\eta(y)\\
& =: & \int \eta(\T(s,y)) d\eta(y).
\end{eqnarray*}
(Recall that $Q$ is the unit square; $Q$ could be replaced by any
bounded convex set containing $E$). Note that a fixed line $\ell$ intersects at most $C 2^n$ rectangles $Q(\xi)$ with $\xi\in\Xi_n$, whence, by using \eqref{eq:ahlfors_regular} and letting $n\rightarrow\infty$, we see that $\eta(\ell)=0$. Therefore
\[
\lim_{t\rightarrow s} \mathbf{1}_{\T(t,y)}(x) =  \mathbf{1}_{\T(s,y)}(x) \quad\mu-\textrm{a.e.},
\]
whence, by the dominated convergence theorem,
\[
\lim_{t\rightarrow s} \eta(\T(t,y)) = \eta(\T(s,y)).
\]
Applying the dominated convergence theorem again we find that the functions $\{\phi_n\}$
are continuous.

It follows from the proof of \cite[Theorem 18.2]{Pesin1997} that
there exists $K_2>0$ such that
\[
K_2^{-1} \tau_n(s) \le \phi_n(s) \le K_2 \tau_n(s).
\]
Therefore we obtain from \eqref{eq:submultiplicative} that
\[
\phi_{m+n}(s) \le K_3 \phi_n(s)\phi_m(\rot^n(s)),
\]
where $K_3 = A\, K_2^3$. Hence, if we let $\overline{\phi}_n = K_3
\phi_n$, we have
\[
\log\ophi_{m+n}(s) \le \log\ophi_n(s) + \log\ophi_m(\rot^n(s)).
\]
We have shown that $\{\log\ophi_n\}$ is a continuous subadditive
cocycle over $([0,\beta),\rot)$. By Theorem \ref{th:furman}, and
taking into account the negative factor $1/\log(a)$, for almost
every $s\in [0,\beta)$ we have
\begin{equation} \label{eq:kingman}
\lim_{n\rightarrow\infty} \frac{\log(\phi_n(s))}{n\log(a)} =
\sup_n \frac{\int \log\phi_n(\zeta)d\leb(\zeta)}{n\log(a)} =:
\widetilde{D}.
\end{equation}
Moreover, for \textit{all} $s\in [0,\beta)$ we have the inequality
\begin{equation} \label{eq:furman}
\liminf_{n\rightarrow\infty} \frac{\log(\phi_n(s))}{n\log(a)} \ge
\widetilde{D}.
\end{equation}

Recall that $\underline{D}(\nu)$ is the supremum of all $\alpha$
such that the $\alpha$-energy $I_\alpha(\nu)$ is finite; see
\eqref{eq:def_energy}. It follows from the potential-theoretic
proof of Marstrand's projection theorem (see e.g. \cite[Chapter
9]{Mattila1995}) that
\[
D(\eta_s) = D(\eta),
\]
for almost every $s\in\mathbb{R}$. Thus $\widetilde{D} = D(\eta) =
\dim_H(E)$.

Lipschitz maps do not increase upper correlation dimension; this
can be easily checked from the definition. Therefore, using
\eqref{eq:dimension_eta},
\[
\overline{D}(\eta_s) \le \overline{D}(\eta) = \dim_H(E),
\]
for \textit{all} $s\in [0,\beta)$. But from \eqref{eq:furman} and
the fact that $\widetilde{D} = \dim_H(E)$ we also get
\[
\underline{D}(\eta_s) \ge \dim_H(E)
\]
for all $s\in [0,\beta)$. Thus for all $s\in [0,\beta)$ we have
$\overline{D}(\eta_s) = \underline{D}(\eta_s) = \dim_H(E)$. This
completes the proof of the theorem. \qed

\section{The case $d_a+d_b>1$} \label{sec:largedimension}

Recall that a real number $\theta$ is a \textbf{Pisot number} if $\theta$ is an algebraic integer, $\theta>1$ and all the algebraic conjugates of $\theta$ have modulus strictly smaller than $1$. The main result of this section is the following:

\begin{theorem} \label{th:example}
Suppose that $0<a<b<1/2$, $\log b/\log a$ is irrational, and $1/a$ and $1/b$ are both Pisot numbers. Then there exists a dense $G_\delta$, and therefore uncountable, set
$B\subset (0,\infty)$, such that if $\lambda\in B$, then
 $\nu_{a,b}^\lambda$ is a singular measure.
\end{theorem}
The proof of  Theorem \ref{th:example} will be given at the end of this section. Let $F(\cdot)$ denote the Fourier transform of a measure,
defined by
\[
F(\mu)(\xi) = \int e^{i x \xi} d\mu(x).
\]
By elementary properties of the Fourier transform, which are still valid for Fourier transforms of measures, we have
\begin{equation} \label{eq:transformconvolution}
F(\nu_{a,b}^\lambda)(\xi) = F(\mu_a)(\xi) \,F(\mu_b)(\lambda \xi).
\end{equation}
Salem proved that $F(\mu_a)(\xi)\nrightarrow 0$ as $\xi\rightarrow \infty$ if and only if $1/a$ is a Pisot number; see e.g. \cite{Salem1963} for a proof, as well as further background on Pisot numbers. Thus, if either $1/a$ or $1/b$ is not Pisot, then
\[
\lim_{\xi\rightarrow\infty} F(\nu_{a,b}^\lambda)(\xi)  =0 \quad\text{ for all } \lambda\in\RR\backslash\{0\}.
\]
In the proof of  Theorem \ref{th:example} we will show a converse of this:
if $1/a$ and $1/b$ are both Pisot (and $0<a,b<1/2$) then, for $\lambda\in B$,
\begin{equation} \label{eq:transformdoesntgotozero}
F(\nu_{a, b}^\lambda)(\xi) \nrightarrow 0 \textrm{ as } \xi
\rightarrow\infty \,.
\end{equation}

Theorem \ref{th:example} provides a counterexample to the principle that
for dynamically defined sets and measures, the set of exceptions to the projection theorems
should be determined by natural algebraic relations.

An open question, due to H. Furstenberg, is whether this principle is
valid for orthogonal projections of simple self-similar sets
like the one-dimensional Sierpi\'{n}ski gasket, defined as
\[
\mathcal{S} = \left\{ \sum_{j=0}^\infty 3^{-i}\omega_i: \omega_i
\in \{ (0,0),(1,0),(0,1)\} \right\}.
\]
The conjecture in this case is that Hausdorff dimension is
preserved for orthogonal projections in directions with irrational
slope. 

By taking $a=1/4$ and $b=1/3$, Theorem \ref{th:example} gives the first example we know for which the
principle described above is known to fail. Here the set of exceptions is
uncountable, and since $\log 4/\log 3$ is irrational, there is no
exact overlap for $\lambda\neq 0$. On the other hand, Theorem
\ref{th:main_result} is one of the few cases in which the
principle has been proved to hold.

By Theorem \ref{th:main_result}, $\nu_{1/3,1/4}^\lambda$ has
correlation, Hausdorff and packing dimension equal to $1$ for all
nonzero $\lambda$. Thus for $\lambda\in  B$ there is a loss of
absolute continuity, but not a dimension drop.
We remark that for  certain
$1$-dimensional measures in $\RR^2$ (Hausdorff measures restricted to purely nonrectifiable
sets of positive finite one-dimensional Hausdorff measure) their projections onto almost every line are
one-dimensional but singular, due to a classical theorem of
Marstrand. The crucial difference is that the measure we are
projecting, $\mu_{1/3}\times \mu_{1/4}$, has dimension strictly
greater than $1$. The same considerations apply to all $0<a,b<1/2$ such that $1/a$ and $1/b$ are both Pisot, $\log b/\log a$ is irrational, and $\dim_H(C_a)+\dim_H(C_b)>1$.

Note also that Theorem \ref{th:example} does not contradict the
result of Senge and Strauss \cite{SengeStraus1971} mentioned in
the introduction, since Senge and Strauss only deal with the case
$\lambda=1$. It is still possible that $\mu_a * \mu_b$ is
absolutely continuous whenever $d_a+d_b>1$ and $\log b/\log
a\notin\mathbb{Q}$, but Theorem \ref{th:example} precludes using
the ideas in the proof of Theorem \ref{th:main_result} to prove
such a result.

For notational reasons it will be convenient to rescale and
translate $\mu_t$ so that the convex hull of its support becomes
$[-1/(1-t),1/(1-t)]$. Since we are rescaling the supports of
$\mu_a$ and $\mu_b$ by a different factor, this change also
induces a rescaling of the parameter $\lambda$, but this does not
affect the statement of Theorem \ref{th:example}.

The advantage of this change of coordinates is that $\mu_t$
becomes the distribution of the random sum
\[
\sum_{j=0}^\infty \pm t^j,
\]
where $\PP(+)=\PP(-)=1/2$ and all the choices are independent.  This in turn yields the well-known expression for the Fourier
transform of $\mu_t$ as an infinite product:
\begin{align*}
F(\mu_t)(\xi) &= \prod_{j=0}^\infty \frac{F(\delta_{t^j})(\xi)+
F(\delta_{-t^j})(\xi)}{2} \\
&= \prod_{j=0}^\infty \frac{e^{ i t^j \xi}+e^{- i t^j \xi}}{2}\\
&= \prod_{j=0}^\infty \cos(t^j \xi),
\end{align*}
where $\delta_x$ denotes the unit Dirac mass at $x$. These
infinite products have been studied intensively, see e.g. \cite{Salem1963}.

The following known lemma describes the set $B$ of parameters $\lambda$
for which we will prove that $\nu_{a,b}^\lambda$ is singular. Since we are not aware of a suitable reference, and the proof is short, we include it for the convenience of the reader.

\begin{lemma} \label{lemma:setofexceptions}
Let $0<a,b<1$ be numbers such that $\log b/\log a$ is irrational, and let $\e>0$ be arbitrary. Let
\[
B = \{ \lambda>0: | \lambda a^{-n} -  b^{-m} | < \e \textrm{ for
infinitely many pairs } (n,m)\in \mathbb{N}^2 \}.
\]
Then $B$ is a dense $G_\delta$ subset of $(0,\infty)$. In
particular, it is uncountable.
\end{lemma}
\begin{proof}
For each $N\in\mathbb{N}$, let
\[
B_N = \{ \lambda>0: | \lambda a^{-n} - b^{-m} | < \e \textrm{ for some
} n, m \ge N \}.
\]
It is clear that $B_N$ is open. We claim that it is also dense in
$(0,\infty)$. Indeed, let $I=(c,d)\subset (0,\infty)$ be any
nonempty interval. Then $I$ meets $B_N$ whenever $a^n b^{-m} \in I$
for some $n,m\ge N$. By taking logarithms this is equivalent to
\[
\frac{\log b}{|\log a|}m - n \in \left(\frac{\log c}{|\log
a|},\frac{\log d}{|\log a|}\right ).
\]
But one can find arbitrarily large integers $n,m$ satisfying this
by the irrationality of $\log b/\log a$.

Noting that $B=\cap_{N\ge 1} B_N$ and applying Baire's Theorem
concludes the proof of the lemma.
\end{proof}

\begin{proof}[Proof of Theorem \ref{th:example}]
We start by recalling some basic facts about Pisot numbers. Let $\theta>2$ be Pisot, and write $\theta_1,\ldots,\theta_r$ for the algebraic conjugates of $\theta$. Since $\theta^n + \sum_{i=1}^r \theta_i^n$ is an integer for all $n\in\mathbb{N}$, we find that
\begin{equation} \label{eq:pisotbasicproperty}
\dist(\theta^n, \mathbb{Z}) \le r \gamma^n,
\end{equation}
for all natural numbers $n$, where
\[
\gamma := \max_{i=1,\ldots,r} |\theta_i| < 1.
\]
Moreover, since $\theta$ is an algebraic integer, so is $\theta^n$ for all $n\in\mathbb{N}$; in particular, $\theta^n$ cannot be of the form $k+1/2$ with $k$ an integer. From these two facts we obtain that the infinite product
\[
F(\mu_{1/\theta})(\pi) = \prod_{j=0}^\infty \cos(\pi \theta^j)
\]
is absolutely convergent. Likewise, since $\theta>2$, then the infinite product $\prod_{j=0}^\infty \cos(\pi \theta^{-j})$ is also absolutely convergent.

The above observations also imply that the number
\begin{equation} \label{eq:definitionepsilon}
\e := \frac{1}{2} \dist\left(\{ b^{-n}:n\in\mathbb{N} \}, \mathbb{Z}+\frac{1}{2}\right)
\end{equation}
is strictly positive.

Now let $B$ be the set given by Lemma
\ref{lemma:setofexceptions}, with $\e$ defined in \eqref{eq:definitionepsilon}. Fix $\lambda\in B$ for the rest of the proof. Using \eqref{eq:transformconvolution}, we see that the Fourier transform of $\nu_{a,b}^\lambda$ is given by
\[
\Phi(\xi) = \prod_{j=0}^\infty \cos(a^j\xi)
\prod_{j=0}^\infty \cos(\lambda b^j\xi) =:
\Phi_1(\xi)\Phi_2(\xi).
\]
Let $N,M\in\mathbb{N}$ be such that $|\lambda a^{-N} -  b^{-M}|<\e$.
We have
\begin{align*}
|\Phi_1(\pi a^{-N})| &= \prod_{j=0}^{\infty} |\cos(\pi a^{-N+j})|\\
&\ge \prod_{j=-\infty}^\infty |\cos(\pi a^{-j})|\\
&=: c_1 > 0,
\end{align*}
by our earlier observations.

Write $\sigma=\lambda a^{-N} - b^{-M}$, and note that by the definition of $\e$, and using that $|\sigma b^j|\le \e$ for $j\ge 0$,
\[
\dist\left( b^{j-M}+\sigma b^j , \mathbb{Z} +\frac{1}{2}\right) \ge \e/2,
\]
for all $j=0,\ldots, M$. Thus, using \eqref{eq:pisotbasicproperty} and that $0<b<1$, we get that the products
\[
\prod_{j=0}^{M-1} |\cos(\pi(b^{j-M}+\sigma b^j) )|
\]
are uniformly bounded below by some constant $c_2>0$ independent of $M$.

Using this, we estimate:
\begin{align*}
|\Phi_2(\pi a^{-N})| &= \prod_{j=0}^\infty |\cos(\pi (b^{-M}+\sigma)b^j)|\\
&=\prod_{j=0}^{M-1} |\cos(\pi(b^{j-M}+\sigma b^j) )|
\prod_{j=M}^\infty
|\cos(\pi b^{j-M}(1+\sigma b^M))| \\
&\ge  c_2 \,|F(\mu_{b})(\pi(1+\sigma b^M))|.
\end{align*}
Recall that $F(\mu_{b})(\pi)\neq 0$; thus
$F(\mu_{b})(\pi(1+\sigma b^M)
)$ is bounded away from zero for
sufficiently large $M$.

We have shown that $\Phi$ is bounded away from zero on the set
\[
\{ \pi a^{-N} : | a^{-N}\lambda -b^{-M}|<\e \textrm{ for some large }
M\in\mathbb{N}\}.
\]
Since, by assumption, the set above is infinite (hence unbounded), this
shows that $\Phi(\xi)\nrightarrow 0$ as $\xi\rightarrow\infty$, so $\nu_{a, b}^\lambda$
 is not absolutely continuous. By The Jessen-Wintner law
of pure types (see, e.g., \cite[Theorem 3.26]{Breiman1968}), it follows that $\nu_{a, b}^\lambda$
is singular, as claimed.
\end{proof}

\section{Remarks and open questions}
\label{sec:conclusion}

We finish the paper with some generalizations, remarks and open
problems.

\textbf{General self-similar sets}. The sets $C_a$ are among the
simplest examples of self-similar sets. Recall that a set
$C\subset\mathbb{R}$ is said to be self-similar if there are
affine maps $f_i(x) = \lambda_i x + t_i$, $i=1,\ldots,m$, with
$|\lambda_i|<1$, such that
\[
C = \bigcup_{i=1}^m f_i(C).
\]
If all the $\lambda_i$ coincide, we say that $C$ is an homogeneous
self-similar set. It is not hard to see that the proof of Theorem
\ref{th:main_result} extends, with minor modifications, to the
restrictions of Hausdorff measure to pairs $C, C'$ of homogeneous
self-similar sets satisfying the strong separation condition. In
the case of sets, it is possible to reduce the general
self-similar case to the homogeneous one, by observing that any
self-similar set contains an homogeneous one of arbitrarily close
dimension; see \cite[Proposition 6]{PeresShmerkin2009}. However,
for measures such reduction does not work, and we do not know if
Theorem \ref{th:main_result} is valid for measures on arbitrary self-similar measures.

\textbf{Bernoulli convolutions}. Notice that the definition of $\mu_t$ as the distribution of a random sum makes sense whenever $t\in(0,1)$; if $t>1/2$ then the support of $\mu_t$ becomes an interval. The family $\mu_t$ for $t\in (0,1)$ is known as the family of \textbf{Bernoulli convolutions}. The proof of Theorem \ref{th:example} applies, with minor modifications, also in the case where $a$ or $b$ are in $(1/2,1)$. The case $a=1/2$ (or $b=1/2$) is exceptional, since $\mu_{1/2}$ is the restriction of Lebesgue measure to its supporting interval. Pisot numbers also play a prominent r\^{o}le in the study of Bernoulli convolutions: the only parameters $t\in (1/2,1)$ for which $\mu_t$ is known to be singular are reciprocal of Pisot numbers; on the other hand, Solomyak proved that $\mu_t$ is absolutely continuous for almost every $t\in (1/2,1)$. The reader is referred to \cite{Solomyak2004} for a proof of these facts and further background on Bernoulli convolutions.

\textbf{The measure of $C_a+\lambda C_b$}. In
\cite{PeresShmerkin2009} it was asked whether $C_a+ C_b$ has positive
Lebesgue measure whenever $d_a+d_b>1$ and $\log b/\log
a\notin\mathbb{Q}$. Theorem \ref{th:example} suggests that $C_a +
\lambda C_b$ may have zero Lebesgue measure for some nonzero
values of $\lambda$ if $1/a$ and $1/b$ are Pisot numbers, but we do not have a proof.
Even in this case, it could still be that $C_a+C_b$ has positive
measure, but one would need a different method of proof.

\textbf{Natural measures on $C_a + C_b$}. Besides $\mu_a
 * \mu_b$, other possible natural measures on $C_a+C_b$ are the restrictions of Hausdorff and packing measures in the appropriate dimension. However, Ero\v{g}lu \cite{Eroglu2007} proved that
\[
\mh^{d_a+d_b}(C_a+C_b) = 0,
\]
for all $a,b$ such that $d_a+d_b\le 1$. We do not know whether
$C_a+C_b$ has positive $d_a+d_b$-dimensional packing measure when
$\log b/\log a$ is irrational (it is easy to see that it is
finite).

\textbf{Sums of more than two central Cantor sets}. Let
$a_1,\ldots, a_k$ be a collection of real numbers in $(0,1/2)$
which is linearly independent over $\mathbb{Q}$. Then
\[
D(\mu_{a_1}*\ldots * \mu_{a_k}) = \min(d_{a_1}+\ldots+d_{a_k},1).
\]
The proof of this is similar to that of Theorem
\ref{th:main_result}. We sketch the main differences. The space
$\XX$ becomes the family of all $k$-tuples of finite words with
elements in $\{0,1\}$. A family $\{ \XX_n\}$ of subsets of $\XX$
is then constructed, with the property that if $\xi\in \XX_n$,
then each parallelepiped $Q(\xi)$, defined in the obvious way, has
size
\[
a_1^n \times a_1^n e^{\rot_2^n(0)/\beta_2} \times\ldots a_1^n
e^{\rot_k^n(0)/\beta_k},
\]
where $\beta_i$ are real numbers (the analogous of $\beta$), and
$\rot_2,\ldots, \rot_n$ are rotations of the circle. It follows
from the hypothesis of linear independence that
\[
\rot(x_2,\ldots,x_k) = (\rot_2(x_2),\ldots,\rot_k(x_k))
\]
is a uniquely ergodic transformation of the $(k-1)$-dimensional
torus. Instead of a single family $\{\YY_n\}$, now $2^{k-1}$
auxiliary families are required (one for each subset of
$\{2,\ldots, k\}$). Details are left to the reader.

\textbf{$L^q$-dimensions}. The $L^q$ dimensions of a measure generalize the concept of correlation dimension. We define them only for $q>1$ since that is the case we are interested in. Let
\[
\mathcal{C}^q_\nu(r) = \int \nu (B(x,r))^{q-1} d\nu(x),
\]
and set
\[
D_q(\nu) = \liminf_{r\rightarrow 0} \frac{\mathcal{C}^q_\nu(r)}{\log(r)}.
\]
Thus $D_2(\nu)$ equals the lower correlation dimension of $\nu$. The $L^q$ dimensions are of fundamental importance in multifractal analysis, see for example \cite{Falconer1997}. In general, there is no projection theorem for $L^q$-dimensions for $q>2$: the $L^q$-dimension can drop for all orthogonal projections, even for very simple measures like arc length on the unit circle. However, the rest of the proof of Theorem \ref{th:main_result} extends to $L^q$ dimensions for any $q>1$. In particular, the analogue of \eqref{eq:submultiplicative} holds for $L^q$-dimensions, with an almost identical proof. Applying Furman's Theorem and the subadditive ergodic theorem as in the proof of Theorem \ref{th:main_result}, we obtain the following result:
\begin{theorem}
Let $q>2$. If
\[
D_q(\mu_a * (\mu_b \circ S_{\lambda_0}^{-1})) > d
\]
for some $\lambda_0>0$, then
\[
D_q(\mu_a * (\mu_b \circ S_\lambda^{-1})) > d
\]
for almost every $\lambda>0$.
\end{theorem}
This is related to the investigations in
\cite{Furstenberg1970}.

\ack{We thank Esa J\"{a}rvenp\"{a}\"{a} for helpful comments}.

\end{document}